\documentclass[11pt]{article}
\usepackage{mathrsfs}
\usepackage{amsmath}
\usepackage{amssymb}
\usepackage{graphicx}
\usepackage{epic}
\renewcommand{\paragraph}{\roman{paragraph}}

\def \la{\lambda}

\newtheorem{theorem}{\scshape \mdseries  Theorem}[section]
\newtheorem{lemma}[theorem]{\scshape \mdseries  Lemma}
\newtheorem{coro}[theorem]{\scshape \mdseries  Corollary}

\begin{document}

\title{\sf  A lower bound for the algebraic connectivity of a graph in terms of the domination number\thanks{
Supported by National Natural Science Foundation of China (11071002, 11371028, 71101002), Program for New Century Excellent
Talents in University (NCET-10-0001), Key Project of Chinese Ministry of Education (210091),
Specialized Research Fund for the Doctoral Program of Higher Education (20103401110002),
Project of Educational Department of Anhui Province (KJ2012B040),
Scientific Research Fund for Fostering Distinguished Young Scholars of Anhui University(KJJQ1001).}
}
\author{Yi-Zheng Fan,\!\!$^{1,}$\thanks{Corresponding author.
 E-mail addresses: fanyz@ahu.edu.cn(Y.-Z. Fan), tansusan1@ahjzu.edu.cn (Y.-Y. Tan).} \ \ Ying-Ying Tan$^2$\\
    {\small  \it $1$. School of Mathematical Sciences, Anhui University, Hefei 230601, P. R. China}\\
  {\small  \it $2$. Department of Mathematics \& Physics, Anhui University of Architecture, Hefei 230601, P. R. China} \\
 }
\date{}
\maketitle

\noindent {\bf Abstract:}
We investigate how the algebraic connectivity of a graph changes by relocating a connected branch from one vertex to another vertex,
and then minimize the algebraic connectivity among all connected graphs of order $n$ with fixed domination number $\gamma \le \frac{n+2}{3}$,
and finally present a lower bound for the algebraic connectivity in terms of the domination number.



\section{Introduction}

Let $G=(V(G),E(G))$ be a simple graph with vertex set $V(G)=\{v_1,v_2,\ldots,v_n\}$ and edge set $E(G)$.
The {\it adjacency matrix} of $G$ is defined to be
a $(0,1)$-matrix $A(G)=[a_{ij}]$,
where $a_{ij}=1$ if $v_i$ is adjacent to $v_j$, and
$a_{ij}=0$ otherwise.
The {\it degree matrix} of $G$ is defined by $D(G)=\hbox{diag}\{d_G(v_1), d_G(v_2),\cdots, d_G(v_n)\}$, where $d_G(v)$ or simply $d(v)$ is the degree
of a vertex $v$ in $G$.
The matrix $L(G)=D(G)-A(G)$ is called the {\it Laplacian matrix} of $G$.
It is known that $L(G)$ is positive semidefinite, and $0$ is the smallest eigenvalue with the all-one vector $\mathbf 1$ as the corresponding eigenvector.
The second smallest eigenvalue of $L(G)$, denoted by $\alpha(G)$, is known as the {\it algebraic connectivity} of $G$ due to Fiedler \cite{fied1}.
The eigenvectors corresponding to $\alpha(G)$, also called {\it Fiedler vectors}, have a nice structural property; see \cite{fied2}.
There are many results on the algebraic connectivity and Fiedler vectors; see e.g. \cite{abr,bappm,fallatk,kirk,merris}.

But, here we consider the perturbations of the algebraic connectivity of a graph under locally changing of the graph.
This has been investigated by Kirkland and Neumnann \cite{kirkn}, Patra and Lal \cite{patl} on trees or weighted trees,
and by Guo \cite{guojm1,guojm2} on general graphs.

Let $G_1$, $G_2$ be two vertex-disjoint graphs, and let $v\in V(G_1)$, $u\in V(G_2)$. The {\it coalescence} of
$G_1$ and $G_2$ with respect to $v$ and $u$, denoted by $G_1(v)\diamond G_2(u)$, is obtained from $G_1$ and
$G_2$ by identifying $v$ with $u$ and forming a new vertex $p$, which is also denoted as $G_1(p)\diamond G_2(p)$.
If a connected graph $G$ can be expressed as $G=G_1(p)\diamond G_2(p)$,
where $G_1$ and $G_2$ are nontrivial subgraphs of $G$ both containing $p$,
then $G_1$ or $G_2$ is called a {\it branch} of $G$ rooted at $p$.
Let $G=G_1(v_2)\diamond G_2(u)$ and $G^*=G_1(v_1)\diamond G_2(u)$, where $v_1$ and $v_2$ are two distinct vertices of $G_1$
and $u$ is a vertex of $G_2$. We say that $G^*$ is obtained from $G$ by {\it relocating $G_2$ from $v_2$ to $v_1$}.

\begin{center}
\includegraphics[scale=.6]{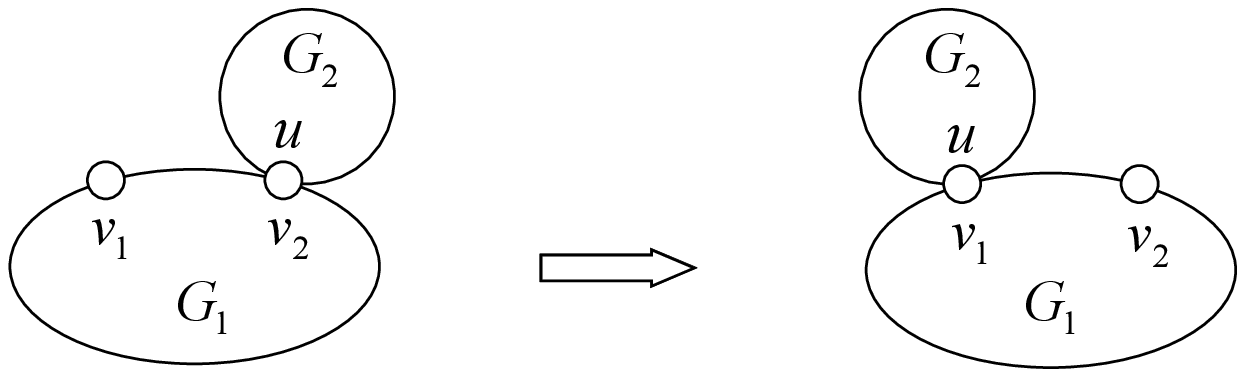}

{\small Fig. 1.1. Relocating $G_2$ from $v_2$ to $v_1$.}
\end{center}

Our problem is: {\it When relocating a branch from one vertex to another vertex, how does the algebraic connectivity change?}
Similar results have been obtained for the least eigenvalue of the adjacency matrix of graphs \cite{fanw} and the signless Laplacian of non-bipartite  graphs \cite{wang}.
In this paper, we first investigate how the algebraic connectivity changes by relocating the branch of a graph,
and then minimize the algebraic connectivity among all connected graphs of order $n$ with fixed domination number $\gamma \le \frac{n+2}{3}$,
and finally present a lower bound for the algebraic connectivity in terms of the domination number.

Recall that a vertex set $S$ of the graph $G$ is called a {\it dominating set} if every vertex of $V(G)\backslash S$ is adjacent to at least one  vertex of $S$.
The {\it domination number} of $G$, denoted by $\gamma(G)$, is the minimum of the cardinalities of all domination sets in $G$.
Lu et al. \cite{lu}, Nikiforov et al. \cite{niki} and Feng \cite{feng} give some upper bounds for the algebraic connectivity of graphs in terms of domination number, respectively.
But no work appears on the lower bound for the algebraic connectivity in terms of domination number.

\section{Perturbation result for the algebraic connectivity}

We first give some preliminary knowledge and notations. 
A graph $G$ is called {\it trivial} if it contains only one vertex;
otherwise, it is called {\it nontrivial}.
Let $G$ be a graph on vertices $v_1,v_2,\ldots, v_n$, and let $x=(x_1, x_2, \ldots, x_n) \in \mathbb{R}^n$.
The vector $x$ can be considered
as a function defined on $V(G)$, which maps each vertex $v_i$ of $G$ to the value $x_i$, i.e. $x(v_i)=x_i$.
If $x$ is an eigenvector of $L(G)$, then it defines on $G$ naturally, i.e. $x(v)$ is the entry of $x$ corresponding to $v$.
One can find
that the quadratic form $x^TL(G)x$ can be written as
$$x^TL(G)x=\sum_{uv \in E(G)}[x(u)-x(v)]^2. \eqno(2.1)$$
The eigenvector equation $L(G)x=\la x$ can be interpreted as
$$[d(v)-\la] x(v)= \sum_{u \in N_G(v)} x(u) \hbox{~ for each~} v \in V(G), \eqno(2.2)$$
where $N_G(v)$ denotes the neighborhood of $v$ in $G$.
In addition, for an arbitrary unit vector $x \in \mathbb{R}^n$ orthogonal to $\mathbf{1}$,
$$\alpha(G)\leq x^TL(G)x, \eqno(2.3)$$
with equality if and only if $x$ is a Fiedler vector of $G$.

The following two lemmas give a nice property of Fiedler vectors of a graph or a tree.

\begin{lemma}{\em \cite{fied2}}\label{fiedler}
Let $G$ be a connected graph with a Fiedler vector $x$.
Let $u$ be a cut vertex of $G$, and let $G_0,G_1,\ldots,G_r\;(r\ge 1)$ be all components of $G-u$.
Then:

{\em (i)} If $x(u)>0$, then exactly one of the components $G_i$ contains a vertex negatively valuated by $x$.
For all vertices $v$ in the remaining components $x(v)>x(u)$.

{\em (ii)} If $x(u)=0$ and there is a component $G_i$ containing both positively and negatively valuated vertices, then
there is exactly one such component, all remaining being zero valuated.

{\em (iii)} If $x(u)=0$ and none component contains  both positively and negatively valuated vertices, then each component $G_i$ contains
either only positively valuated, or negatively valuated, or only zero valuated vertices.

\end{lemma}

\begin{lemma}{\em \cite{fied2}} \label{fiedler2}
Let $T$ be a tree with a Fiedler vector $x$.
Then exactly one of the two cases occurs:

Case A. All values of $x$ are nonzero. Then $T$ contains exactly one edge $pq$ such that $x(p)>0$ and $x(q)<0$.
The values in vertices along any path in $T$ which starts in $p$ and does not contain $q$ strictly increase, the values in vertices
along any path starting in $q$ and not containing $p$ strictly decrease.

Case B. The set $N_0=\{v: x(v)=0\}$ is non-empty. Then the graph induced by $N_0$ is connected and there is exactly one vertex $z \in N_0$ having
at least one neighbor not belonging to $N_0$. The values along any path in $T$ starting in $z$ are strictly increasing, or strictly decreasing, or zero.
\end{lemma}

If the Case B in Lemma \ref{fiedler2} occurs, the vertex $z$ is called the {\it characteristic vertex}, and $T$ is called a {\it Type I} tree;
otherwise, $T$ is called a {\it Type II} tree in which case the edge $pq$ is called the {\it characteristic edge}.
The characteristic vertex or characteristic edge of a tree is independent of the choice of Fiedler vectors; see \cite{merris}.

\begin{lemma}{\em \cite{kirk}}\label{cutvertex}
 Let $G$ be a connected graph with a cut vertex $v$. Then $\alpha(G) \le 1$, with equality if and
only if $v$ is adjacent to all other vertices of $G$.
\end{lemma}

Now we give a perturbation result on the algebraic connectivity of a graph by relocating one branch from one vertex to another vertex.

\begin{lemma}\label{relocate}
Let $G_1$ be a connected graph containing at least two vertices $v_1$, $v_2$, and let $G_2$ be a nontrivial connected graph containing a
vertex $u$. Let $G=G_1(v_2)\diamond G_2(u)$ and $G^*=G_1(v_1)\diamond G_2(u)$. If there exist a Fiedler vector $x$ of $G$ such that
$x(v_1)\geq x(v_2) \ge 0$ and all vertices in $G_2$ are nonnegatively valuated by $x$, then
$$\alpha(G^*) \leq \alpha(G),$$
with equality if and only if $x(v_1)= x(v_2)=0$, $\sum_{w\in N_{G_2}(u)}x(w)=0$, and $x$ is also a Fiedler vector of $G^*$.
\end{lemma}

{\it Proof:}
Assume that $x$ has unit length.
Suppose $G_1$ has $n_1$ vertices and $G_2$ has $n_2$ vertices.
Let $n:=n_1+n_2-1$, the number of vertices of $G$.
Let $y$ be a vector defined on the graph $G^*$ such that $y(w)=x(w)+[x(v_1)-x(v_2)]$ for each $w \in V(G_2)\backslash\{u\}$, and $y(w) =x(w)$ for all remaining vertices $w$.
Then $\mathbf{1}^T y=(n_2-1)[x(v_1)-x(v_2)]$.
Define $$z=y-\frac{(n_2-1)[x(v_1)-x(v_2)]}{n}\mathbf{1}.$$
Then $\mathbf{1}^T z=0$, $z^TL(G^*)z=x^TL(G)x$, and
$$z^Tz=1+2[x(v_1)-x(v_2)]\sum_{w \in V(G_2)\backslash\{u\}}x(w)+\frac{n_1(n_2-1)}{n}[x(v_1)-x(v_2)]^2 \ge 1.$$
So we have $\alpha(G^*) \leq \alpha(G)$.

If $\alpha(G^*) = \alpha(G)$, then $z^Tz=1$, which implies that $x(v_1)=x(v_2)$.
So $z=x$, which is also a Fiedler vector of $G^*$.
If $x(v_2)>0$, then by Lemma \ref{fiedler}(i) the component containing negatively valuated vertices is contained in $G_1-u$ as $x(v_1)=x(v_2)$.
So for each vertices $w \in V(G_2)\backslash\{u\}$, $x(w) > x(v_2)$ also by Lemma \ref{fiedler}(i).
However, if considering the eigenvector equation (2.2) of $G$ and $G^*$ on the vertex $v_1$, we will have
$$d_{G_2}(u)x(v_1)=\sum_{w\in N_{G_2}(u)}x(w)>d_{G_2}(u)x(v_2)=d_{G_2}(u)x(v_1),$$
a contradiction.
So, $x(v_1)=x(v_2)=0$ and $\sum_{w\in N_{G_2}(u)}x(w)=0$.
The sufficiency for  $\alpha(G^*) = \alpha(G)$ is easily verified as $x^TL(G^*)x=x^TL(G)x$ and $x$ is a  Fiedler vector of $G^*$.
\hfill $\blacksquare$

{\bf Remark:}
The result in Lemma \ref{relocate} does not hold if $x(v_1),x(v_2)$ have different signs.
For example, let $G:=T(k,l,2)$ be the graph as listed in Fig. 3.1, where $k \ge 1$ and $l \ge 1$.
By Lemma \ref{cutvertex}, $\alpha(G)<1$.
If letting $x$ be a Fiedler vector of $G$, then $x$ contains no zero entries by eigenvector equation (2.2).
So, by Lemma \ref{fiedler2}, $x(v_1)x(v_2)<0$, and all pendant vertices attached at $v_1$ (respectively, $v_2$) have the same sign as $v_1$ (respectively, $v_2$).
Relocating the pendant star attached at $v_2$ to $v_1$, or relocating the pendant star attached at $v_1$ to $v_2$,
we always get a star $G^*$ such that $\alpha(G^*)=1>\alpha(G)$.

\section{Minimizing the algebraic connectivity}
For convenience a graph is called {\it minimizing} among a certain class of graphs, if its algebraic connectivity attains the minimum 
among all graphs in such class.
In this section we first characterize the minimizing tree(s) among all trees of order $n$ with domination number $\gamma$, and then
extend the result to general graphs, where $n \ge 3\gamma-2$.


Denote by $T(k,l,d)\;(k \ge l)$ the tree of order $n$ obtained from a path $P_d$ by attaching $k$ and $l$ pendant edges respectively at
its two endpoints, where $n=k+l+d$; see Fig. 2.1.
If $k=0$ or $l=0$, then no pendant edges are attached at $v_1$ or $v_d$.

If $d \ge 2$, $k \ge 1$ and $l \ge 1$, then $\alpha(T(k,l,d))<1$ by Lemma \ref{cutvertex}.
By eigenvector equation and by Lemma \ref{fiedler2}, the star attached at $v_1$ has all vertices nonzero valuated with the same sign,
and the star attached at $v_d$ has all vertices nonzero valuated with the same sign but different to the above star.

\begin{center}
\includegraphics[scale=.5]{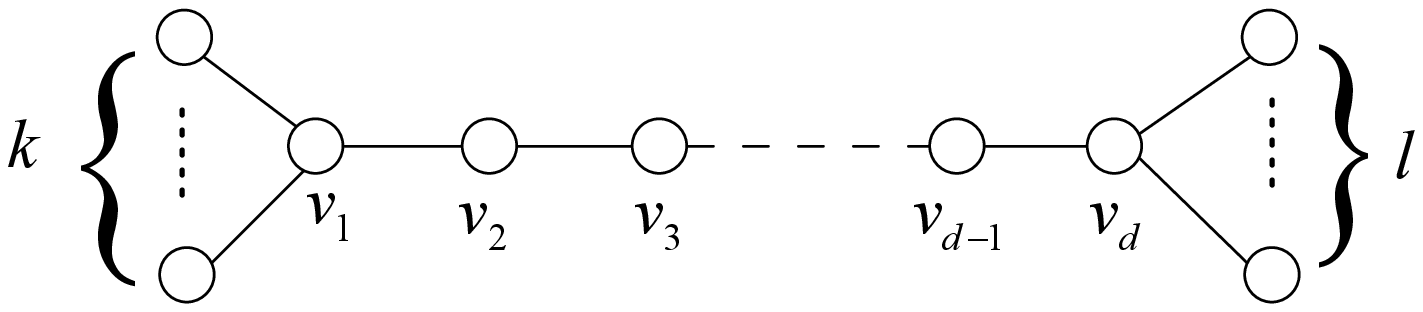}

{\small Fig. 2.1. The tree $T(k,l,d)$}
\end{center}

\begin{lemma}{\em \cite{fallatk}} \label{treedia}
Among all trees of order $n$ and diameter $d+1$, the tree $T_d:=T(\lceil \frac{n-d}{2} \rceil, \lfloor \frac{n-d}{2} \rfloor, d)$ is the unique
graph with minimum algebraic connectivity.
\end{lemma}

\begin{lemma}\label{ineq1}
{\em (1)}  If $k \ge 2$, $\alpha(T(k,l,d)) > \alpha(T(k-1,l,d+1))$;

\noindent {\em (2)} if $l \ge 2$  $\alpha(T(k,l,d)) > \alpha(T(k,l+1,d+1))$.
\end{lemma}

{\it Proof:}
If $d=1$, then $T(k,l,d)$ is a star and both inequalities hold by Lemma \ref{cutvertex}.
Assume $d \ge 2$.
Let $x$ be a Fiedler vector of $T(k,l,d)$.
By Lemma \ref{fiedler2} and the fact $\alpha(T(k,l,d))<1$, without loss of generality, assume that the star $S$ attached at $v_1$ is positively valuated by $x$.
Let $u_1,u_2,\ldots,u_k$ be all pendant vertices adjacent to $v_1$, and let $S'$ be the sub-star of $S$ consisting of vertices $u_2,\ldots,u_k$ and $v_1$.
Then $x(u_1) > x(v_1)>0$ by the eigenvector equation (2.2).
Relocating $S'$ from $v_1$ to $u_1$, we will arrive at a graph isomorphic to $T(k-1,l,d+1)$.
By Lemma \ref{relocate}, $\alpha(T(k,l,d)) > \alpha(T(k-1,l,d+1))$.
The second inequality can be proved similarly. \hfill $\blacksquare$

\begin{lemma}\label{ineq2}
Let $\gamma(T_d))=\gamma \ge 2$, where $n \ge d+2$ and $n \ge 3\gamma+1$.
Then $ 3\gamma-4 \le d \le 3\gamma -2$.
Furthermore,
$$\alpha(T_d))
\ge \alpha(T_{3\gamma -2})),$$
with equality if and only if $d=3\gamma -2$.
\end{lemma}

{\it Proof:}
Clearly $\gamma(T_d)=\lceil \frac{d+2}{3} \rceil$.
So  $\gamma-1 < \frac{d+2}{3} \le \gamma$, and thus $ 3\gamma-4 \le d \le 3\gamma -2$.
If $d=3\gamma-3$, noting $n \ge 3\gamma+1$, then by Lemma \ref{ineq1} and Lemma \ref{treedia},
$$
\alpha(T_d)  >  \alpha(T(\lceil (n-d-2)/2 \rceil, \lfloor (n-d)/2 \rfloor, d+1))\ge   \alpha(T_{3\gamma -2}).
$$
If $d=3\gamma-4 $, also by Lemma \ref{ineq1},
\begin{align*}
\alpha(T_d) & > \alpha(T(\lceil (n-d-2)/2 \rceil, \lfloor (n-d)/2 \rfloor, 3\gamma-3))\\
& > \alpha(T(\lceil (n-d-2)/2 \rceil, \lfloor (n-d-2)/2 \rfloor, 3\gamma-2)\\
&= \alpha(T_{3\gamma -2}).
\end{align*}
\hfill $\blacksquare$

\begin{lemma}\label{ineq3}
If $\gamma_1 < \gamma_2$ and $n \ge 3\gamma_2+1$,
then $\alpha(T_{3\gamma_1-2}) >  \alpha(T_{3\gamma_2-2})$.
\end{lemma}

{\it Proof:}
By Lemma \ref{ineq1} and Lemma \ref{treedia},
\begin{align*}
\alpha(T_{3\gamma_1-2}) & > \alpha(T(\lceil (n-(3\gamma_1-2)-2)/2 \rceil, \lfloor (n-(3\gamma_1-2))/2 \rfloor, 3\gamma_1-1))\\
& > \alpha(T(\lceil (n-(3\gamma_1-2)-2)/2 \rceil, \lfloor (n-(3\gamma_1-2)-2)/2 \rfloor, 3\gamma_1))\\
& >\alpha(T(\lceil (n-(3\gamma_1-2)-2-2)/2 \rceil, \lfloor (n-(3\gamma_1-2)-2)/2 \rfloor, 3\gamma_1+1))\\
& \ge \alpha(T(\lceil (n-(3\gamma_1+1))/2 \rceil, \lfloor (n-(3\gamma_1+1))/2 \rfloor, 3\gamma_1+1))\\
&=\alpha(T_{3(\gamma_1+1)-2}).
\end{align*}
The result follows by induction on the domination number.\hfill $\blacksquare$

\begin{theorem}\label{maintree}
Among all trees of order $n$ and domination number $\gamma$, where $n \ge 3\gamma-2$,
 the tree $T_{3\gamma-2}$ is the unique
graph with minimum algebraic connectivity.
\end{theorem}

{\it Proof:}
The result clear holds for $\gamma=1$.
In addition, if $n \in \{3\gamma-2,3\gamma-1,3\gamma\}$, then $T_{3\gamma-2}$ will one of $T(0,0,3\gamma-2)$, $T(1,0,3\gamma-2)$ and $T(1,1,3\gamma-2)$.
Surely $T_{3\gamma-2}=P_n$, and the result holds as $P_n$ is the unique minimizing graph among all connected graphs of order $n$.
Suppose $\gamma \ge 2$ and $n \ge 3\gamma+1$.
Let $T$ be a minimizing tree. 
A {\it pendent star} of $T$ is a maximal subtree of $T$ induced on pendant vertices together with the quasi-pendent vertex to which they all are attached.
If $T$ has exactly two pendant stars, then $T=T(k,l,d)$ for some $k,l,d$, where $d \ge 2$.
By Lemma \ref{treedia}, $k=\lceil \frac{n-d}{2} \rceil$ and $l=\lfloor \frac{n-d}{2} \rfloor$.
The result follows by Lemma \ref{ineq2}.

Now suppose that $T:=T_0$ has more than two pendant stars, which has $p_0$ pendent vertices and $q_0$ quasi-pendent vertices.
Let $x$ be a Fiedler vector of $T_0$.
If $T_0$ is of Type I, then there exist at least one zero pendant star $S$ attached at some vertex say $u$,
  and at least one positive quasi-pendant vertex $w$.
Relocating the zero star $S$ at $u$ to $w$, we will arrive at a new tree $T_1$ such that $\alpha(T_1) < \alpha(T_0) $ by Lemma \ref{relocate}.
Note that $\gamma(T_1) \le \gamma(T_0)$. In fact, $\gamma(T_1) < \gamma(T_0)$; otherwise we will get a contradiction to the fact that $T_0$ is minimizing.
If $T$ is of Type II, then there exist at least two pendant stars $S_1,S_2$ both being positive or negative valuated by $x$, attached at $u_1,u_2$ respectively.
Without loss of generality, assume $S_1,S_2$ are both positive and $x(u_1) \ge x(u_2) >0$.
Relocating $S_2$ from $u_2$ to $u_1$, we also arrive at a new tree $T_1$ such that $\alpha(T_1) < \alpha(T_0) $ by Lemma \ref{relocate} and $\gamma(T_1) < \gamma(T_0)$.

Repeat the above procession on $T_1$ if $T_1$ has more than two pendant stars and continue a similar discussion to the resulting trees.
Note that from the $k$-th step to the $(k+1)$-th step, either $p_{k+1}=p_k$ and $q_{k+1}=q_k-1$, or $p_{k+1}=p_k+1$ and $q_{k+1}=q_k$.
So the above procession will be terminated at the $n$-th step in which the tree $T_n$ has exactly two pendant stars, i.e. $T_n=T(k,l,d)$ for some $k,l,d$, where $d \ge 2$.
Hence
$$\alpha(T)=\alpha(T_0)>\alpha(T_1)> \cdots >\alpha(T_n), ~~\gamma(T)=\gamma(T_0)>\gamma(T_1)>\cdots >\gamma(T_n).$$
Noting that $\gamma(T_{3\gamma-2})=\gamma$ and $\gamma(T_d)=\gamma(T_n)$, by Lemma \ref{treedia} and Lemma \ref{ineq2}, we have
$$\alpha(T_{3\gamma-2}) \ge \alpha(T) > \alpha(T_n) \ge \alpha(T_d) \ge \alpha(T_{3\gamma(T_n)-2}).$$
However, since $\gamma(T_n)<\gamma(T)=\gamma$, by Lemma \ref{ineq3}, we have $\alpha(T_{3\gamma-2}) <  \alpha(T_{3\gamma(T_n)-2})$, a contradiction.
So this case cannot happen and the result follows.\hfill $\blacksquare$

\begin{lemma} \label{spantree}
Let $G$ be a connected graph of order $n$ and domination number $\gamma$. Then $G$ contains a spanning tree with domination number $\gamma$.
\end{lemma}

{\it Proof:}
If $\gamma=1$, the result holds obviously.
Now suppose $\gamma \ge 2$.
Let $U=\{u_1,u_2,\ldots,u_\gamma\}$ be a dominating set of $G$ of size $\gamma$, let $W=V(G) \backslash U$.
Let $B$ be a bipartite spanning subgraph of $G$, which is obtained by deleting all possible edges within $U$ or $W$.

First assume that $B$ is connected. Then there exist two vertices in $U$, say $u_1$ and $u_2$, such that $N_B(u_1)\cap N_B(u_2)\neq \emptyset$.
Assume that $w_1 \in N_B(u_1)\cap N_B(u_2)$.
Deleting all edges between $u_2$ and the vertices of $(N_B(u_1)\cap N_B(u_2))\backslash \{w_1\}$ (if it is nonempty),
we will get a subgraph $B_1$ of $B$ such that $u_2$ shares exactly one neighbor with $u_1$.
If $U\backslash \{u_1,u_2\} \ne \emptyset$,
noting that $B_1$ is also connected, there exists one vertex $w_2 \in N_B(u_1)\cap N_B(u_2)$ such that $w_2$ is adjacent to one vertex, say $u_3$ in $U\backslash \{u_1,u_2\}$.
Deleting all edges between $u_3$ and  the vertices of $(N_{B_1}(u_3) \backslash  \{w_2\}) \cap (N_B(u_1)\cup N_B(u_2))$,
we will get a subgraph $B_2$ of $B_1$ such that $u_3$ shares exactly one neighbor with $u_1$ or $u_2$.
Repeating the above process, we will arrive at a subgraph $B_{\gamma-1}$ of $B$ such that for each $i=2,3, \ldots,\gamma$,
$u_i$ shares exactly one neighbor with $u_1$, $u_2$, $\ldots$, or $u_{i-1}$.
So $B_{\gamma-1}$ is a tree with domination number $\gamma$, as desired.

Next suppose that $B$ is not connected.
Let $B_1,B_2,\ldots,B_k$ be the components of $B$ with bipartitions  $(U_1, W_1), (U_2, W_2), \ldots, (U_k, W_k)$ respectively.
By the above discussion, each $B_i$ contains a spanning tree $T_i$ such that $\gamma(T_i)=\gamma(B_i)$ for $i=1,2,\ldots,k$.
Since $G$ is connected, there exists a spanning tree $T$ of $G$ obtained from $T_i$'s by adding $k-1$ edges between $U_i$s and $U_j$s, or $W_i$s and $W_j$s.
Surely $\gamma(T)=\gamma$ and the result follows.
\hfill $\blacksquare$

\begin{lemma}{\em \cite{shao}} \label{shao}
Let $G$ be a connected graph with a pendant star on vertices $v_0,v_1, \ldots, v_s \;(s \ge 2)$, where $v_0$ is the vertex to which the star is attached.
Let $G'$ be obtained from $G$ by arbitrarily adding edges among $v_1,\ldots,v_s$.
If $\alpha(G) \ne 1$, then $\alpha(G)=\alpha(G')$.
\end{lemma}

\begin{lemma}{\em \cite{merris2}} \label{merris2}
Let $G_1,G_2$ be two graphs of order $r,s$ respectively.
If the eigenvalues of $L(G_1)$ are $0,\mu_1,\mu_2,\ldots,\mu_{r-1}$ and the
eigenvalues of $L(G_2)$ are $0,\nu_1,\nu_2,\ldots,\nu_{s-1}$, then the eigenvalues of $G_1 \vee G_2$ are
$0, \mu_1+s,\mu_2+s,\ldots,\mu_{r-1}+s, \nu_1+r,\nu_2+r,\ldots,\nu_{s-1}+r, r+s$, where $G_1 \vee G_2$ denotes the
graph obtained from the union of $G_1$ and $G_2$ by adding all possible edges between every vertex of $G_1$ and
each of $G_2$.
\end{lemma}

\begin{theorem}\label{main}
Among all connected graphs of order $n$ and domination number $\gamma$, where $n \ge 3\gamma-2$,
if $G$ is a minimizing graph, then $G$ is one of the following graph:

\noindent
{\em (1)} if $\gamma =1$, then $G$ contains a cut vertex that is adjacent to all other vertices.

\noindent
{\em (2)} if  $\gamma \ge 2$, then $G$ is obtained from the tree $T_{3\gamma-2}$ by arbitrarily adding some (or none) edges among the pendant vertices in a same pendant star.
\end{theorem}

{\it Proof:}
Let $G$ be a minimizing graph.
First suppose that $\gamma \ge 2$.
If $n \in \{3\gamma-2,3\gamma-1,3\gamma\}$, Then the result obviously holds as the path $P_n$ is the unique minimizing graph.
Now suppose $n \ge 3\gamma+1$.
By Lemma \ref{spantree}, $G$ contains a spanning tree $T$ also with domination number $\gamma$.
By Theorem \ref{maintree},
$$\alpha(G) \ge \alpha(T) \ge \alpha(T_{3\gamma-2}).$$
If $\alpha(G) =\alpha(T_{3\gamma-2})$, then $T=T_{3\gamma-2}$ also by Theorem \ref{maintree}.

Returning to the origin graph $G$, which is obtained from $T_{3\gamma-2}$ possibly by adding some edges.
Assume that $E(G)\backslash E(T_{3\gamma-2})\neq \emptyset$.
Let $x$ be a unit Fiedler vector of $G$. Then
\begin{align*}
\alpha(G) &= \sum_{uv \in E(G)}[x(u)-x(v)]^2\\
& =\sum_{uv \in E(T_{3\gamma-2})}[x(u)-x(v)]^2+\sum_{uv \in E(G)\backslash E(T_{3\gamma-2})}[x(u)-x(v)]^2\\
& \geq \sum_{uv \in E(T_{3\gamma-2})}[x(u)-x(v)]^2 \geq \alpha(T_{3\gamma-2}).
\end{align*}
Since $\alpha(G) =\alpha(T_{3\gamma-2})$, $x$ is also  a Fiedler vector of $T_{3\gamma-2}$, and $x(u)-x(v)=0$ for each edge
$uv \in E(G)\backslash E(T_{3\gamma-2})$.
By Lemma \ref{fiedler2}, $u,v$ are both the pendent vertices lying in a same pendant star.
So the necessity follows.
The sufficiency follows from the result Lemma \ref{shao}.

If $\gamma=1$, then $G$ contains a spanning star centered at $u$, i.e. $G=\{u\} \vee H$ for some graph $H$.
By Lemma \ref{merris2}, $\alpha(G) \ge 1$ with equality if and only if $H$ has the zero eigenvalue with multiplicity at least two, i.e. $H$ is disconnected.
So $u$ is a cut vertex of $G$.
\hfill $\blacksquare$

\begin{lemma}{\em \cite{kirkn}} \label{bound}
Suppose that $d \ge 3$, $ k \ge 1$, $l \ge 1$ and $n:=k+l+d-1$.
Then
$$\alpha(T(k,l,d-1)) \ge \left(\frac{nd}{4}-\frac{2n+d^2-4d-5}{8}\right)^{-1}.$$
\end{lemma}

\begin{coro}
Let $G$ be a connected graph of order $n$ and domination number $\gamma \le \frac{n+2}{3}$.
Then $$\alpha(G) \ge \frac{8}{6n\gamma+18\gamma-9\gamma^2-4n}.$$
\end{coro}

{\it Proof:}
By Theorem \ref{main}, $\alpha(G) \ge \alpha(T_{3\gamma-2})$.
If $\gamma=1$, surely $\alpha(T_{3\gamma-2})=1 > \frac{8}{6n\gamma+18\gamma-9\gamma^2-4n}.$
If $\gamma \ge 2$ and $n \ge 3\gamma$, the result follows if taking $d=3\gamma-1$ in Lemma \ref{bound}.
If $n$ equals $3\gamma-2$ or  $3\gamma-1$, noting that in this case $T_{3\gamma-2}=P_n$ and
$\alpha(P_{3\gamma-2}) > \alpha(P_{3\gamma-1}) > \alpha(P_{3\gamma})$.
So the result also holds. \hfill $\blacksquare$

\end{document}